\documentclass[11pt,a4paper]{article}

\usepackage{inputenc}
\usepackage{amsmath}
\usepackage{bm}
\usepackage{bbold}
\usepackage{amsthm}
\usepackage{enumerate}

\usepackage{hyperref}
\usepackage{breakurl}

\title{Tropical optimization techniques in multi-criteria decision making with Analytical Hierarchy Process\thanks{UKSim-AMSS 11th European Modelling Symposium on Computer Modelling and Simulation (EMS 2017) / Ed. by D.~Al-Dabass, Z.~Xie, A.~Orsoni, A.~Pantelous. IEEE, 2017. pp.~38-43.}}

\author{N. Krivulin\thanks{Saint Petersburg State University, Faculty of Mathematics and Mechanics, 28 Universitetsky Ave., St.~Petersburg, 198504, Russia, nkk@math.spbu.ru.}
\and
S. Sergeev\thanks{University of Birmingham, School of Mathematics, Edgbaston B15 2TT, UK, s.sergeev@bham.ac.uk}
}

\date{}

\newtheorem{theorem}{Theorem}
\newtheorem{lemma}[theorem]{Lemma}
\newtheorem{corollary}[theorem]{Corollary}

\theoremstyle{definition}
\newtheorem{example}{Example}

\begin{document}

\maketitle

\begin{abstract}
We apply methods and techniques of tropical optimization to develop a new theoretical and computational framework for the implementation of the Analytic Hierarchy Process in multi-criteria problems of rating alternatives from pairwise comparison data. The framework involves the Chebyshev approximation of pairwise comparison matrices by consistent matrices in the logarithmic scale. We reduce the log-Chebyshev approximation to multidimensional tropical optimization problems, and offer complete direct solutions to the problems in the framework of tropical mathematics. The results obtained provide a closed-form solution to the rating problem of interest as either a unique score vector (up to a positive factor) or as a set of different score vectors. To handle the problem when the solution is not unique, we develop tropical optimization techniques to find those vectors from the solution set that are the most and least differentiating between the alternatives with the highest and lowest scores, and thus can be well representative of the entire solution.
\\

\\

\end{abstract}

\section{Introduction}

Tropical (idempotent) mathematics, which deals with the theory and applications of algebraic systems with idempotent operations \cite{Golan2003Semirings,Heidergott2006Maxplus,Mceneaney2006Maxplus,Itenberg2007Tropical,Gondran2008Graphs,Butkovic2010Maxlinear,Maclagan2015Introduction}, is widely used as a coherent analytical framework to solve problems in engineering, operations research and computer science. Tropical optimization presents an important research domain in this area, focused on optimization problems that are formulated and solved in the tropical mathematics setting.

Methods and techniques of tropical optimization are applied to solve many well-known and new optimization problems in various fields, including decision making \cite{Elsner2004Maxalgebra,Elsner2010Maxalgebra,Gursoy2013Theanalytic,Tran2013Pairwise,Gavalec2015Decision}. Specifically, in \cite{Gursoy2013Theanalytic}, the Analytic Hierarchy Process (AHP) method in multi-criteria problems of rating alternatives from pairwise comparison data is investigated from the tropical point of view.

The traditional AHP method \cite{Saaty1977Scaling,Saaty1980Analytic,Saaty2013Onthemeasurement} consists of two principal levels of pairwise comparisons: the upper level, where the relative importance of criteria is estimated, and the lower level, where the relative quality of choices is evaluated with respect to each criterion. The final decision is made by combining the rates of all choices computed on the lower level and the weights of all criteria on the higher level.

The rates of choices with respect to each criterion are found through the rank-one approximation of pairwise comparison matrices, typically by using principal (Perron) eigenvector methods \cite{Saaty1977Scaling,Saaty1980Analytic,Saaty1984Comparison,Saaty2013Onthemeasurement}, and, sometimes, by other techniques, including least squares or logarithmic least squares methods \cite{Saaty1984Comparison,Chu1998Ontheoptimal,Saaty1984Comparison,Barzilai1997Deriving,Farkas2003Consistency,Gonzalezpachon2003Transitive}. The weights of criteria can be evaluated in the same manner from pairwise comparisons or obtained in a different way.

More specifically, assume that there are $m$ criteria and $n$ choices. Given a pairwise comparison matrix $\bm{A}_{k}=(a_{ij}^{(k)})$ of order $n$, and the weight $w_{k}$ for each criterion $k$, the vector of the priorities of all choices $\bm{x}=(x_{i})$ is calculated as
\begin{equation*}
\bm{x}
=
\sum_{k=1}^{m}w_{k}\bm{x}_{k}, 
\end{equation*}
where $\bm{x}_{k}$ is the vector of rates, obtained from $\bm{A}_{k}$.

In this paper, we develop a new theoretical and computational framework for the implementation of AHP, based on the tropical optimization techniques proposed in \cite{Krivulin2015Rating,Krivulin2016Using,Krivulin2018Methods}. The new AHP method, which we offer and investigate below, aims to find a rank-one matrix that should be the closest, in the sense of the maximum of weighted log-Chebyshev distances, to all the pairwise comparison matrices corresponding to different criteria. In this method, the vector of priorities $\bm{x}=(x_{i})$ is a solution of the following optimization problem
\begin{equation*}
\min_{\bm{x}}\ 
\max_{i,j=1}^{n}\max_{k=1}^{m}(w_{k}a_{ij}^{(k)})x_{j}/x_{i}.
\end{equation*}

The solution obtained as a result of the tropical AHP method is, in general, non-unique. However, inconsistency and ambiguity of determining the pairwise preferences are inherent to the process of forming the pairwise comparison matrices, and therefore it seems quite natural that the method ends up with a set of solution vectors rather than with just one vector as in the 
ordinary AHP.

To make the non-unique result tractable and useful for practice, we focus on two kinds of solutions that can be considered, in some sense, as the best and worst solutions. The solution set is characterized by vectors that are the most and least differentiating between the choices with the highest and lowest priorities. We obtain the most and least differentiating vectors by solving the span (range) seminorm maximization and minimization problems
\begin{equation*}
\max_{\bm{x}}\
(\max_{i=1}^{n}x_{i})/\max_{j=1}^{n}(1/x_{j})),
\quad
\min_{\bm{x}}\
(\max_{i=1}^{n}x_{i})(\max_{j=1}^{n}(1/x_{j})).
\end{equation*}

We formulate the above optimization problems in the framework of tropical mathematics as tropical optimization problems, and then offer direct complete solutions to these problems, represented in a compact vector form. To illustrate the results obtained, we demonstrate the numerical solution of a multi-criteria decision problem from \cite{Saaty1977Scaling}.

\section{Log-Chebyshev Approximation Based AHP}

In this section, we describe a new approach to develop an AHP decision scheme that is based on the rank-one log-Chebyshev approximation of pairwise comparison matrices. 

\subsection{Log-Chebyshev Approximation of Comparison Matrices}

Consider the problem of evaluating the rates of $n$ choices from the pairwise comparison of the choices. The outcome of the comparison is described by a square symmetrically reciprocal matrix $\bm{A}=(a_{ij})$, where $a_{ij}$ specifies the relative priority of choice $i$ over $j$, and satisfies the condition $a_{ij}=1/a_{ji}>0$ for all $i,j$. The pairwise comparison matrix $\bm{A}$ is called consistent if its entries are transitive, that is, if they satisfy the equality $a_{ij}=a_{ik}a_{kj}$ for all $i,j,k$.

For each consistent matrix $\bm{A}$, there is a positive vector $\bm{x}=(x_{i})$ whose elements completely determine the entries of $\bm{A}$ by the relation $a_{ij}=x_{i}/x_{j}$, which, in particular, means that $\bm{A}$ is a matrix of unit rank. Provided that the matrix $\bm{A}$ is consistent, its corresponding vector $\bm{x}$, which can be readily obtained from $\bm{A}$, directly represents, up to a positive factor, the individual preferences of choices in question.

Since the pairwise comparison matrices, encountered in practice, are generally inconsistent, the solution usually involves approximating these matrices by consistent matrices. The approximation with the principal (Perron) eigenvector as well as the least squares or the logarithmic least squares approximation are often used as solution approaches.

Consider another approach, which is based on the approximation of a pairwise comparison matrix $\bm{A}=(a_{ij})$ by a consistent matrix $\bm{X}=(x_{ij})$ in the log-Chebyshev sense, where the approximation error is measured with the Chebyshev metric on the logarithmic scale. Taking into account that both matrices $\bm{A}$ and $\bm{X}$ have positive entries and that the logarithmic function (to the base more than one) is monotonically increasing, the error can be written as
\begin{equation*}
\max_{i,j=1}^{n}|\log a_{ij}-\log x_{ij}|
=
\log\max_{i,j=1}^{n}\max\{a_{ij}/x_{ij},x_{ij}/a_{ij}\}.
\end{equation*}

Observing that the minimization of the logarithm is equivalent to the minimization of its argument and that $a_{ij}=1/a_{ji}$ and $x_{ij}=x_{i}/x_{j}$, we replace the last logarithm by 
$\max_{i,j=1}^{n}\max\{a_{ij}/x_{ij},x_{ij}/a_{ij}\}=\max_{i,j=1}^{n}a_{ij}x_{j}/x_{i}$. This reduces the approximation problem under study to the solution, with respect to the unknown vector of priorities $\bm{x}=(x_{i})$, of the optimization problem 
\begin{equation}
\min_{\bm{x}}\ (\max_{i,j=1}^{n}a_{ij}x_{j}/x_{i}).
\label{P-minxaijxjxi}
\end{equation}

\subsection{Weighted Approximation Under Several Criteria}

Suppose the priorities of choices are evaluated based on pairwise comparisons according to $m$ criteria, each having a given weight. For each criterion $k$, we denote the pairwise comparison matrix by $\bm{A}_{k}=(a_{ij}^{(k)})$ and the positive weight by $w_{k}$. To determine the priority vector $\bm{x}=(x_{i})$, we minimize the maximum of the functions $\max_{i,j=1}^{n}a_{ij}^{(k)}x_{j}/x_{i}$, taken with the weights $w_{k}$ for all $k$, which is given by
\begin{equation*}
\max_{k=1}^{m}w_{k}(\max_{i,j=1}^{n}a_{ij}^{(k)}x_{j}/x_{i})
=
\max_{i,j=1}^{n}\max_{k=1}^{m}(w_{k}a_{ij}^{(k)})x_{j}/x_{i}.
\end{equation*}

Then, we introduce the matrix $\bm{C}=(c_{ij})$ with the entries
\begin{equation}
c_{ij}
=
\max_{k=1}^{m}w_{k}a_{ij}^{(k)},
\label{E-cij-maxwkaijk}
\end{equation}
and arrive at the problem in the form of \eqref{P-minxaijxjxi}. Note that the solution of problem \eqref{P-minxaijxjxi} with the matrix $\bm{C}$ can be considered as a modification of the basic AHP scheme, in which the log-Chebyshev approximation is used instead of the  principal eigenvector method, and the weights of criteria are incorporated into the lower level evaluation of choices.

\subsection{Most and Least Differentiating Priority Vectors}

The priority vectors that solve problem \eqref{P-minxaijxjxi} are in general non-unique up to multiplication by a positive factor. Below, we offer an approach to reduce the non-unique solution to two representative vectors, which can be considered, in some sense, as the best and worst solutions.

Assume that problem \eqref{P-minxaijxjxi} has a set $S$ of solutions $\bm{x}=(x_{i})$ rather than a unique one. Since the main purpose of evaluating priorities is to differentiate choices, we find those solutions that are the most and least differentiating between the choices with the highest and lowest priorities. The calculation of the most and least differentiating vectors involves determining the boundaries within which the contrast ratio
\begin{equation*}
(\max_{i=1}^{n}x_{i})/(\min_{j=1}^{n}x_{j})
=
(\max_{i=1}^{n}x_{i})(\max_{j=1}^{n}(1/x_{j}))
\end{equation*}
lies to find the vectors $\bm{x}$, which solve the problem of span (range) seminorm maximization
\begin{equation}
\max_{\bm{x}\in S}\ (\max_{i=1}^{n}x_{i})(\max_{j=1}^{n}(1/x_{j})),
\label{P-maxxmaxximax1xj}
\end{equation}
and the problem of the seminorm minimization
\begin{equation}
\min_{\bm{x}\in S}\ (\max_{i=1}^{n}x_{i})(\max_{j=1}^{n}(1/x_{j})).
\label{P-minxmaxximax1xj}
\end{equation}

In subsequent sections, we represent problems \eqref{P-minxaijxjxi}, \eqref{P-maxxmaxximax1xj} and \eqref{P-minxmaxximax1xj} in terms of tropical mathematics, and give direct and explicit solutions, which are ready for immediate computation.

\section{Elements of tropical mathematics}

We start with a brief overview of basic definitions and notation of tropical (idempotent) algebra to provide a formal framework for describing tropical optimization techniques, used below to develop a tropical implementation of AHP. Further details on tropical mathematics can be found, e.g., in the recent works \cite{Golan2003Semirings,Heidergott2006Maxplus,Mceneaney2006Maxplus,Itenberg2007Tropical,Gondran2008Graphs,Butkovic2010Maxlinear,Maclagan2015Introduction}. 

Consider the system $(\mathbb{R}_{+},\oplus,\otimes,0,1)$, where $\mathbb{R}_{+}$ denotes the set of non-negative reals, which is closed under two operations: addition $\oplus$ defined as maximum with neutral element $0$, and multiplication $\otimes$ defined as usual with neutral element $1$. Addition is idempotent as $x\oplus x=\max(x,x)=x$ for all $x\in\mathbb{R}_{+}$. Multiplication is distributive over addition, and invertible since each $x\ne0$ has an inverse $x^{-1}$ such that $x\otimes x^{-1}=xx^{-1}=1$. The system under consideration is called the idempotent semifield or the max-algebra and denoted by $\mathbb{R}_{\max}$. In what follows, the multiplication sign $\otimes$ is omitted. The power notation has the usual meaning.

The set of matrices over $\mathbb{R}_{+}$ with $m$ rows and $n$ columns is denoted by $\mathbb{R}_{+}^{m\times n}$. A matrix with all zero entries is the zero matrix. Tropical matrix operations are defined by the conventional formulae where the scalar operations $\oplus$ and $\otimes$ play the role of the usual addition and multiplication.

The multiplicative conjugate transpose of a nonzero matrix $\bm{A}=(a_{ij})\in\mathbb{R}_{+}^{m\times n}$ is the matrix $\bm{A}^{-}=(a_{ij}^{-})\in\mathbb{R}_{+}^{n\times m}$ with the entries $a_{ij}^{-}=a_{ji}^{-1}$ if $a_{ji}\ne0$, and $a_{ij}^{-}=0$ otherwise.

Column vectors with $n$ entries from $\mathbb{R}_{+}$ form the set $\mathbb{R}_{+}^{n}$. The vector with all elements equal to $0$ is the zero vector $\bm{0}$. The vector with all elements equal to $1$ is denoted by $\bm{1}$.

A vector $\bm{b}\in\mathbb{R}_{+}^{n}$ is linearly dependent on the system of vectors $\bm{a}_{1},\ldots,\bm{a}_{n}\in\mathbb{R}_{+}^{n}$ if $\bm{b}=x_{1}\bm{a}_{1}\oplus\cdots\oplus x_{n}\bm{a}_{n}$ for some numbers $x_{1},\ldots,x_{n}\in\mathbb{R}_{+}$. In particular, a vector $\bm{b}$ is collinear with a vector $\bm{a}$, if $\bm{b}=x\bm{a}$ for a number $x$.

The multiplicative conjugate transpose of a nonzero column vector $\bm{x}=(x_{i})$ is the row vector $\bm{x}^{-}=(x_{i}^{-})$, where $x_{i}^{-}=x_{i}^{-1}$ if $x_{i}\ne0$, and $x_{i}^{-}=0$ otherwise.  

A matrix $\bm{A}$ is of rank $1$ if and only if $\bm{A}=\bm{x}\bm{y}^{T}$, where $\bm{x}$ and $\bm{y}$ are positive column vectors. A matrix $\bm{A}$ that satisfies the condition $\bm{A}^{-}=\bm{A}$ is called symmetrically reciprocal (or reciprocal). A reciprocal matrix $\bm{A}$ is of rank $1$ if and only if $\bm{A}=\bm{x}\bm{x}^{-}$, where $\bm{x}$ is a positive column vector.

Consider square matrices in the set $\mathbb{R}_{+}^{n\times n}$. For any matrix $\bm{A}$ and integer $p>0$, the tropical (or $\max$-algebraic) power notation is routinely defined by the inductive rule $\bm{A}^{p}=\bm{A}^{p-1}\bm{A}$, where $\bm{A}^{0}=\bm{I}$ is the usual identity matrix. 

The tropical ($\max$-algebraic) spectral radius of a matrix $\bm{A}=(a_{ij})$ is computed as the maximum cycle geometric mean of the matrix entries, which is given by
\begin{equation}
\lambda
=
\bigoplus_{1\leq k\leq n}\bigoplus_{1\leq i_{1},\ldots,i_{k}\leq n}(a_{i_{1}i_{2}}a_{i_{2}i_{3}}\cdots a_{i_{k}i_{1}})^{1/k}.
\label{E-lambda-ai1i2ai2i3aiki1}
\end{equation}

Provided that $\lambda\leq1$, the asterate operator (the Kleene star) maps the matrix $\bm{A}$ onto the matrix
\begin{equation*}
\bm{A}^{\ast}
=
\bm{I}\oplus\bm{A}\oplus\cdots\oplus\bm{A}^{n-1}.
\end{equation*}

\section{Tropical Optimization Problems}

In this section, we consider tropical optimization problems, which provide the basis for our tropical implementation of AHP. First, assume that, given a matrix $\bm{A}\in\mathbb{R}_{+}^{n\times n}$, we need to find vectors $\bm{x}\in\mathbb{R}_{+}^{n}$ that solve the problem
\begin{equation}
\min_{\bm{x}}\ 
\bm{x}^{-}\bm{A}\bm{x}.
\label{P-minxxAx}
\end{equation}

A complete, direct solution to the problem is obtained in \cite{Krivulin2015Extremal} as follows (see, also \cite{Butkovic2010Maxlinear} and references therein).
\begin{lemma}
\label{L-minxxAx}
Let $\bm{A}$ be a matrix with tropical spectral radius $\lambda>0$. Then, the minimum value in \eqref{P-minxxAx} is equal to $\lambda$, and all positive solutions are given by
\begin{equation*}
\bm{x}
=
(\lambda^{-1}\bm{A})^{\ast}\bm{u},
\quad
\bm{u}
>
\bm{0}.
\end{equation*}
\end{lemma}

We now suppose that $\bm{A}_{1},\ldots,\bm{A}_{m}\in\mathbb{R}_{+}^{n\times n}$ are given matrices, and $w_{1},\ldots,w_{m}\in\mathbb{R}_{+}$ are given numbers. The problem is to find vectors $\bm{x}\in\mathbb{R}_{+}^{n}$ that attain the minimum
\begin{equation}
\min_{\bm{x}}\ 
\bigoplus_{k=1}^{m}w_{k}\bm{x}^{-}\bm{A}_{k}\bm{x}.
\label{P-minxwkxAkx}
\end{equation}

As a direct consequence of the previous result, we have the following solution \cite{Krivulin2016Using}.
\begin{corollary}
\label{C-minxwkxAkx}
Let $\bm{A}_{1},\ldots,\bm{A}_{m}$ be matrices and $w_{1},\ldots,w_{m}$ be positive numbers such that the matrix $\bm{C}=\bigoplus_{k=1}^{m}w_{k}\bm{A}_{k}$ has the tropical spectral radius $\mu>0$.

Then, the minimum value in \eqref{P-minxwkxAkx} is equal to $\mu$, and all positive solutions are given by
\begin{equation*}
\bm{x}
=
(\mu^{-1}\bm{C})^{\ast}\bm{u},
\quad
\bm{u}
>
\bm{0}.
\end{equation*}
\end{corollary}

We conclude this section with the solution of two optimization problems, which are formulated to maximize and minimize a common objective function. Let $\bm{B}\in\mathbb{R}_{+}^{n\times l}$ be a given matrix, and $\bm{u}\in\mathbb{R}_{+}^{l}$ be an unknown vector. We start with the maximization problem
\begin{equation}
\max_{\bm{u}}\ 
\bm{1}^{T}\bm{B}\bm{u}(\bm{B}\bm{u})^{-}\bm{1}.
\label{P-maxu1BuBu1}
\end{equation}

The next result obtained in \cite{Krivulin2018Methods} (see also \cite{Krivulin2016Maximization}) offers a complete solution to the problem.
\begin{lemma}\label{L-maxu1BuBu1}
Let $\bm{B}$ be a positive matrix and $\bm{B}_{lk}$ be the matrix derived from $\bm{B}=(\bm{b}_{j})$ by fixing the entry $b_{lk}$ and replacing the others by $0$.

Then, the maximum in \eqref{P-maxu1BuBu1} is equal to $\Delta=\bm{1}^{T}\bm{B}\bm{B}^{-}\bm{1}$, and all positive solutions are given by
\begin{equation*}
\bm{u}
=
(\bm{I}\oplus\bm{B}_{lk}^{-}\bm{B})\bm{v},
\quad
\bm{v}
>
\bm{0},
\end{equation*}
where $k$ and $l$ satisfy the condition $\bm{1}^{T}\bm{b}_{k}b_{lk}^{-1}=\Delta$.
\end{lemma}

Finally, we consider the minimization problem in the form
\begin{equation}
\min_{\bm{u}}\ 
\bm{1}^{T}\bm{B}\bm{u}(\bm{B}\bm{u})^{-}\bm{1},
\label{P-minu1BuBu1}
\end{equation}
which has a solution described as follows \cite{Krivulin2017Algebraic,Krivulin2018Methods}.
\begin{lemma}\label{L-minu1BuBu1}
Let $\bm{B}$ be a matrix without zero rows and columns, and let $\widehat{\bm{B}}$ be the sparsified matrix derived from $\bm{B}=(\bm{b}_{j})$ by setting to $0$ each entry $b_{ij}$ which is below the threshold $\Delta^{-1}\bm{1}^{T}\bm{b}_{j}$, where $\Delta=(\bm{B}(\bm{1}^{T}\bm{B})^{-})^{-}\bm{1}$.

Let $\mathcal{B}$ be the set of matrices obtained from $\widehat{\bm{B}}$ by fixing one nonzero entry in each row and setting the others to $0$.

Then, the minimum in \eqref{P-minu1BuBu1} is equal to $\Delta$, and all positive solutions are given by
\begin{equation*}
\bm{u}
=
(\bm{I}\oplus\Delta^{-1}\bm{B}_{1}^{-}\bm{1}\bm{1}^{T}\bm{B})\bm{v},
\quad
\bm{v}>\bm{0},
\quad
\bm{B}_{1}
\in
\mathcal{B}.
\end{equation*}
\end{lemma}

\section{Applications to AHP}

We are now in a position to describe our tropical implementation of AHP. Let us consider a multi-criteria decision problem to rate $n$ alternatives (choices) from pairwise comparisons with respect to $m$ criteria. Suppose that $\bm{A}_{0}\in\mathbb{R}_{+}^{m\times m}$ is a matrix of pairwise comparisons of the criteria, and $\bm{A}_{1},\ldots,\bm{A}_{m}\in\mathbb{R}_{+}^{n\times n}$ are matrices of pairwise comparisons of the alternatives for every criterion. Given the above matrices, the problem consists in finding vectors $\bm{x}\in\mathbb{R}_{+}^{n}$ of scores (rates, priorities) of the alternatives.   

At first, we need to evaluate the relative importance of the criteria by solving problem \eqref{P-minxaijxjxi} with the matrix $\bm{A}_{0}$. In terms of the tropical semifield $\mathbb{R}_{\max}$, problem \eqref{P-minxaijxjxi} takes the form of \eqref{P-minxxAx}. Therefore, we can apply Lemma~\ref{L-minxxAx} to obtain the vector of weights $\bm{w}=(\lambda^{-1}\bm{A}_{0})^{\ast}\bm{v}$, where $\lambda$ is the tropical spectral radius of the matrix $\bm{A}_{0}$, and $\bm{v}$ is any positive vector.

The next step is the evaluation of the priorities of alternatives, which involves the solution of problem \eqref{P-minxaijxjxi} with the matrix $\bm{C}$ defined by \eqref{E-cij-maxwkaijk}. After translation into the language of tropical algebra, we have problem \eqref{P-minxwkxAkx}. Corollary~\ref{C-minxwkxAkx} offers the solution to the problem in the form $\bm{x}=(\mu^{-1}\bm{C})^{\ast}\bm{u}$, where $\bm{C}=\bigoplus_{k=1}^{m}w_{k}\bm{A}_{k}$, $\mu$ is the spectral radius of the matrix $\bm{C}$, and $\bm{u}$ is any positive vector.

If the obtained solution $\bm{x}=\bm{B}\bm{u}$, where $\bm{B}=(\mu^{-1}\bm{C})^{\ast}$ is not unique (up to a positive factor), we need to solve problems \eqref{P-maxxmaxximax1xj} and \eqref{P-minxmaxximax1xj} to determine the most and least differentiating priority vectors. In terms of tropical algebra, the objective function in these problems is represented as $\bm{1}^{T}\bm{x}\bm{x}^{-}\bm{1}=\bm{1}^{T}\bm{B}\bm{u}(\bm{B}\bm{u})^{-}\bm{1}$. In this case, problem \eqref{P-maxxmaxximax1xj} turns into problem \eqref{P-maxu1BuBu1}, which has the solution provided by Lemma~\ref{L-maxu1BuBu1}. Similarly, problem \eqref{P-minxmaxximax1xj} becomes problem \eqref{P-minu1BuBu1} with the solution given by Lemma~\ref{L-minu1BuBu1}.
 
We conclude this section with an example intended to illustrate the computational technique involved in the tropical implementation of AHP, described above. 

\begin{example}
Consider an example \cite{Saaty1977Scaling} where a plan for vacation is to be selected. The places considered include: (1) short trips from Philadelphia (i.e., New York, Washington, Atlantic City, New Hope, etc.), (2) Quebec, (3) Denver, (4) California. The problem is to evaluate the places with respect to the following criteria: (1) cost of the trip from Philadelphia, (2) sight-seeing opportunities, (3) entertainment (doing things), (4) way of travel, (5) eating places. 

The comparison matrix of criteria for vacation is given by
\begin{equation*}
\bm{A}_{0}
=
\begin{pmatrix}
1 & 1/5 & 1/5 & 1 & 1/3
\\
5 & 1 & 1/5 & 1/5 & 1
\\
5 & 5 & 1 & 1/5 & 1
\\
1 & 5 & 5 & 1 & 5
\\
3 & 1 & 1 & 1/5 & 1
\end{pmatrix}.
\end{equation*}

The pairwise comparison matrices of vacation sites with respect to the criteria are:
\begin{gather*}
\bm{A}_{1}
=
\begin{pmatrix}
1 & 3 & 7 & 9
\\
1/3 & 1 & 6 & 7
\\
1/7 & 1/6 & 1 & 3
\\
1/9 & 1/7 & 1/3 & 1
\end{pmatrix},
\qquad
\bm{A}_{2}
=
\begin{pmatrix}
1 & 1/5 & 1/6 & 1/4
\\
5 & 1 & 2 & 4
\\
6 & 1/2 & 1 & 6
\\
4 & 1/4 & 1/6 & 1
\end{pmatrix},
\\
\bm{A}_{3}
=
\begin{pmatrix}
1 & 7 & 7 & 1/2
\\
1/7 & 1 & 1 & 1/7
\\
1/7 & 1 & 1 & 1/7
\\
2 & 7 & 7 & 1
\end{pmatrix},
\qquad
\bm{A}_{4}
=
\begin{pmatrix}
1 & 4 & 1/4 & 1/3
\\
1/4 & 1 & 1/2 & 3
\\
4 & 2 & 1 & 3
\\
3 & 1/3 & 1/3 & 1
\end{pmatrix},
\\
\bm{A}_{5}
=
\begin{pmatrix}
1 & 1 & 7 & 4
\\
1 & 1 & 6 & 3
\\
1/7 & 1/6 & 1 & 1/4
\\
1/4 & 1/3 & 4 & 1
\end{pmatrix}.
\end{gather*}

To solve the problem, we first evaluate the weights of criteria by applying Lemma~\ref{L-minxxAx}. We take the matrix $\bm{A}_{0}$ and find its spectral radius. Using the formula at \eqref{E-lambda-ai1i2ai2i3aiki1} yields
\begin{equation*}
\lambda
=
(a_{14}a_{43}a_{32}a_{21})^{1/4}
=
5^{3/4}.
\end{equation*}

Furthermore, we consider the matrix
\begin{equation*}
\lambda^{-1}\bm{A}_{0}
=
\begin{pmatrix}
5^{-3/4} & 5^{-7/4} & 5^{-7/4} & 5^{-3/4} & 5^{-3/4}/3
\\
5^{1/4} & 5^{-3/4} & 5^{-7/4} & 5^{-7/4} & 5^{-3/4}
\\
5^{1/4} & 5^{1/4} & 5^{-3/4} & 5^{-7/4} & 5^{-3/4}
\\
5^{-3/4} & 5^{1/4} & 5^{1/4} & 5^{-3/4} & 5^{1/4}
\\
3\cdot5^{-3/4} & 5^{-3/4} & 5^{-3/4} & 5^{-7/4} & 5^{-3/4}
\end{pmatrix},
\end{equation*}
and calculate its powers to obtain the matrix
\begin{multline*}
(\lambda^{-1}\bm{A}_{0})^{\ast}
=
\bm{I}\oplus\lambda^{-1}\bm{A}_{0}\oplus\lambda^{-2}\bm{A}_{0}^{2}\oplus\lambda^{-3}\bm{A}_{0}^{3}\oplus\lambda^{-4}\bm{A}_{0}^{4}
\\
=
\begin{pmatrix}
1 & 5^{-1/4} & 5^{-1/2} & 5^{-3/4} & 5^{-1/2}
\\
5^{1/4} & 1 & 5^{-1/4} & 5^{-1/2} & 5^{-1/4}
\\
5^{1/2} & 5^{1/4} & 1 & 5^{-1/4} & 1
\\
5^{3/4} & 5^{1/2} & 5^{1/4} & 1 & 5^{1/4}
\\
3\cdot5^{-3/4} & 3\cdot5^{-1} & 3\cdot5^{-5/4} & 3\cdot5^{-3/2} & 3\cdot5^{-5/4}
\end{pmatrix}.
\end{multline*}

Since all columns of the matrix obtained are collinear, any column of this matrix can serve as the weight vector. We take the first column, and use its elements as coefficients to combine the matrices $\bm{A}_{1},\ldots,\bm{A}_{5}$ into one matrix
\begin{multline*}
\bm{C}
=
\bm{A}_{1}
\oplus
5^{1/4}\bm{A}_{2}
\oplus
5^{1/2}\bm{A}_{3}
\oplus
5^{3/4}\bm{A}_{4}
\oplus
(3\cdot5^{-3/4})\bm{A}_{5}
\\
=
\begin{pmatrix}
5^{3/4} & 7\cdot5^{1/2} & 7\cdot5^{1/2} & 9
\\
5^{5/4} & 5^{3/4} & 6 & 3\cdot5^{3/4}
\\
4\cdot5^{3/4} & 2\cdot5^{3/4} & 5^{3/4} & 3\cdot5^{3/4}
\\
3\cdot5^{3/4} & 7\cdot5^{1/2} & 7\cdot5^{1/2} & 5^{3/4}
\end{pmatrix}.
\end{multline*}

We now apply Corollary~\ref{C-minxwkxAkx}. The evaluation of the tropical spectral radius of the matrix $\bm{C}$ yields
\begin{equation*}
\mu
=
(c_{13}c_{31})^{1/2}
=
2\cdot5^{5/8}7^{1/2}.
\end{equation*}

After calculating the matrix $\mu^{-1}\bm{C}$ and its powers, we combine them to construct the matrix
\begin{equation*}
(\mu^{-1}\bm{C})^{\ast}
=
\bm{I}\oplus\mu^{-1}\bm{C}\oplus\mu^{-2}\bm{C}^{2}\oplus\mu^{-3}\bm{C}^{3}
=
\begin{pmatrix}
1 & r/4 & r/4 & 3/4
\\
3/r & 1 & 3/4 & 3/r
\\
4/r & 1 & 1 & 3/r
\\
1 & r/4 & r/4 & 1
\end{pmatrix},
\end{equation*}
where we denote $r=2\cdot7^{1/2}5^{-1/8}\approx4.3272>4$

Next, we eliminate those columns of the matrix obtained, which are linearly dependent on others. Since the first column is collinear with the third, one of them, say the third column, can be removed. Observing that none of the other columns is linearly dependent on others, we form the matrix 
\begin{equation*}
\bm{B}
=
\begin{pmatrix}
1 & r/4 & 3/4
\\
3/r & 1 & 3/r
\\
4/r & 1 & 3/r
\\
1 & r/4 & 1
\end{pmatrix}
\end{equation*}
to represent a complete solution as the set of priority vectors
\begin{equation*}
\bm{x}
=
\bm{B}\bm{u},
\quad
\bm{u}>\bm{0}.
\end{equation*}

To find the most differentiating solution, we need to apply Lemma~\ref{L-maxu1BuBu1}. We start with the calculation
\begin{gather*}
\bm{1}^{T}\bm{b}_{1}
=
1,
\quad
\bm{1}^{T}\bm{b}_{2}
=
r/4,
\quad
\bm{1}^{T}\bm{b}_{3}
=
1,
\\
\bm{b}_{1}^{-}\bm{1}
=
r/3,
\quad
\bm{b}_{2}^{-}\bm{1}
=
1,
\quad
\bm{b}_{3}^{-}\bm{1}
=
r/3,
\end{gather*}
and then obtain $\Delta=\bm{1}^{T}\bm{B}\bm{B}^{-}\bm{1}=r/3$.

The condition $\bm{1}^{T}\bm{b}_{k}b_{lk}^{-1}=\Delta$ holds if we take either $k=1$ and $l=2$, $k=3$ and $l=2$, or $k=3$ and $l=3$. 

First, assume that $k=1$ and $l=2$. We form the matrices
\begin{gather*}
\bm{B}_{21}
=
\begin{pmatrix}
0 & 0 & 0
\\
3/r & 0 & 0
\\
0 & 0 & 0
\\
0 & 0 & 0
\end{pmatrix},
\qquad
\bm{B}_{21}^{-}\bm{B}
=
\begin{pmatrix}
1 & r/3 & 1
\\
0 & 0 & 0
\\
0 & 0 & 0
\end{pmatrix},
\\
\bm{B}(\bm{I}\oplus\bm{B}_{21}^{-}\bm{B})
=
\begin{pmatrix}
1 & r/3 & 1
\\
3/r & 1 & 3/r
\\
4/r & 4/3 & 4/r
\\
1 & r/3 & 1
\end{pmatrix}.
\end{gather*}

Observing that all columns in the last matrix are collinear to each other, we take one of them, say the first, to write one of the most differentiating solutions as
\begin{equation*}
\bm{x}
=
\begin{pmatrix}
1
\\
3/r
\\
4/r
\\
1
\end{pmatrix}
u,
\quad
u>0,
\quad
r
=
2\cdot7^{1/2}5^{-1/8}.
\end{equation*}

Specifically, by setting $u=1$, we arrive at the priority vector $\bm{x}\approx(1.0000, 0.6933, 0.9244, 1.0000)^{T}$. This vector specifies the order of choices as $(4)\equiv(1)\succ(3)\succ(2)$.

Next, we examine the case where $k=3$ and $l=2$. In a similar way, we obtain the vector 
\begin{equation*}
\bm{x}
=
\begin{pmatrix}
3/4
\\
3/r
\\
3/r
\\
1
\end{pmatrix}
u,
\quad
u>0,
\quad
r
=
2\cdot7^{1/2}5^{-1/8},
\end{equation*}
which suggests another most differentiating solution.

If $u=1$, then $\bm{x}\approx(0.7500, 0.6933, 0.6933, 1.0000)^{T}$, which puts the choices in the order $(4)\succ(1)\succ(3)\equiv(2)$.

It is not difficult to verify that the case with $k=3$ and $l=3$ introduces no other solutions than those already found.

We now turn to an application of Lemma~\ref{L-minu1BuBu1} to derive the least differentiating vector of priorities. First, we calculate
\begin{equation*}
\bm{B}(\bm{1}^{T}\bm{B})^{-}
=
\begin{pmatrix}
1
\\
4/r
\\
4/r
\\
1
\end{pmatrix},
\quad
\Delta
=
(\bm{B}(\bm{1}^{T}\bm{B})^{-})^{-}\bm{1}
=
r/4.
\end{equation*}

The threshold and sparsified matrices are given by
\begin{equation*}
\Delta^{-1}\bm{1}\bm{1}^{T}\bm{B}
=
\begin{pmatrix}
4/r & 1 & 4/r
\\
4/r & 1 & 4/r
\\
4/r & 1 & 4/r
\\
4/r & 1 & 4/r
\end{pmatrix},
\qquad
\widehat{\bm{B}}
=
\begin{pmatrix}
1 & r/4 & 0
\\
0 & 1 & 0
\\
4/r & 1 & 0
\\
1 & r/4 & 1
\end{pmatrix}.
\end{equation*}

We have to examine the matrices in $\mathcal{B}$ obtained from $\widehat{\bm{B}}$ by leaving one nonzero entry in each row. Consider the matrix
\begin{equation*}
\bm{B}_{1}
=
\begin{pmatrix}
1 & 0 & 0
\\
0 & 1 & 0
\\
4/r & 0 & 0
\\
1 & 0 & 0
\end{pmatrix},
\end{equation*}
and successively obtain
\begin{gather*}
\bm{B}_{1}^{-}\bm{1}
=
\begin{pmatrix}
r/4
\\
1
\\
0
\end{pmatrix},
\qquad
\bm{B}_{1}^{-}\bm{1}\bm{1}^{T}
=
\begin{pmatrix}
r/4 & r/4 & r/4
\\
1 & 1 & 1
\\
0 & 0 & 0
\end{pmatrix},
\\
\bm{I}\oplus\Delta^{-1}\bm{B}_{1}^{-}\bm{1}\bm{1}^{T}
=
\begin{pmatrix}
1 & 1 & 1
\\
4/r & 1 & 4/r
\\
0 & 0 & 1
\end{pmatrix}.
\end{gather*}

It remains to calculate the matrix
\begin{equation*}
\bm{B}(\bm{I}\oplus\Delta^{-1}\bm{B}_{1}^{-}\bm{1}\bm{1}^{T})
=
\begin{pmatrix}
1 & r/4 & 1
\\
4/r & 1 & 4/r
\\
4/r & 1 & 4/r
\\
1 & r/4 & 1
\end{pmatrix}.
\end{equation*}

Since all columns of the last matrix are collinear, we take one of them to write the least differentiating solution as
\begin{equation*}
\bm{x}
=
\begin{pmatrix}
1
\\
4/r
\\
4/r
\\
1
\end{pmatrix}
u,
\quad
u>0,
\quad
r
=
2\cdot7^{1/2}5^{-1/8}.
\end{equation*}

Under the condition that $u=1$, we have the priority vector $\bm{x}\approx(1.0000, 0.9244, 0.9244, 1.0000)^{T}$. The vector obtained arranges the choices as $(4)\equiv(1)\succ(3)\equiv(2)$.

Calculation with the other matrices in the set $\mathcal{B}$ gives the same results, and thus is omitted.

As one can see, all solutions indicate the highest score of the fourth choice. The score assigned to the first choice is the same or lower. The third choice has the same or higher score, than the second choice, and both of them always have a lower score than the first. Combining all solutions yields the order of choices in the form $(4)\succeq(1)\succ(3)\succeq(2)$.

Finally, note that the results obtained above with the tropical modification of AHP are quite different from those offered by the classical AHP method. Specifically, the order of choices, found in \cite{Saaty1977Scaling}, is given by $(1)\succ(3)\succ(4)\succ(2)$. 
\end{example}

\section{Conclusion}
We have described a new tropical implementation of AHP, in which the resulting vector of priorities appears as a solution to a problem of tropical optimization. Such solution can be non-unique, but this seems to be quite natural, since the ambiguity of determining the pairwise preferences is inherent to the process of forming the pairwise comparison matrices. If the solution is not unique, we apply tropical optimization techniques to choose particular solutions that are the most and least differentiating between the choices with the highest and lowest scores. The resulting order of choices is in general different from that obtained with the conventional AHP as illustrated with an example from \cite{Saaty1977Scaling}.

The new tropical implementation of AHP provides opportunities for further research. Specifically, better algorithms for finding the most and least differentiating solutions are desirable, as well as some comparison between this new implementation of AHP and the existing ones.

\section*{Acknowledgment}
This work was supported in part by the EPSRC grant EP/P019676/1, and by the Russian Foundation for Humanities grant 16-02-00059.

\bibliographystyle{abbrvurl}

\bibliography{Tropical_optimization_techniques_in_multi-criteria_decision_making_with_Analytical_Hierarchy_Process}

\end{document}